\documentclass[12pt]{article}

\usepackage{amsmath,amssymb,amsfonts,theorem,makeidx,latexsym,epsfig,subfigure}

\usepackage{color}

\newtheorem{defn}{Definition}[section]

\newtheorem{lemma}[defn]{Lemma}

{\theorembodyfont{\rmfamily}

\newtheorem{ex}[defn]{Example}}

\newtheorem{thm}[defn]{Theorem}

\newtheorem{prop}[defn]{Proposition}

\newtheorem{cor}[defn]{Corollary}

\numberwithin{equation}{section}

\newcommand{\h}{{\cal H}}

\newcommand{\ltr}{ L^2(\mathbb R) }

\newcommand{\mn}{\mathbb N}

\newcommand{\mr}{\mathbb R}

\newcommand{\mz}{\mathbb Z}

\newcommand{\mc}{\mathbb C}

\newcommand{\mts}{ \{E_{mb}T_{na}g \}_{m,n \in \mz}}

\def\bp{{\noindent\bf Proof. \ }}

\def\ep{\hfill$\square$\par\bigskip}

\def\bqs{\begin{equation}}

\def\eqs{\tag*{$\square$}\end{equation}\par\bigskip}

\def\la{\langle}

\def\ra{\rangle}

\def\ga{\gamma}

\def\supp{\text{supp}}

\def\bop{\begin{op}\rm}

\def\eop{\end{op}}

\def\bee{\begin{eqnarray}}

\def\ene{\end{eqnarray}}

\def\bes{\begin{eqnarray*}}

\def\ens{\end{eqnarray*}}

\def\bei{\begin{itemize}}

\def\eni{\end{itemize}}

\def\bt{\begin{thm}}

\def\et{\end{thm}}

\def\bc{\begin{cor}}

\def\ec{\end{cor}}

\def\bpr{\begin{prop}}

\def\epr{\end{prop}}

\def\bl{\begin{lemma}}

\def\el{\end{lemma}}

\def\bd{\begin{defn}}

\def\ed{\end{defn}}

\def\bex{\begin{ex}}

\def\enx{\end{ex}}

\def\bfi{\begin{fig}}

\def\efi{\end{fig}}

\def\inr{\int_{-\infty}^\infty}

\title{B-spline approximations of the Gaussian, their Gabor frame properties, and
approximately dual frames
}

\date{\today}

\author{Ole Christensen, Hong Oh Kim, Rae Young Kim}

\begin{document}

\maketitle

\begin{abstract}
We prove that Gabor systems generated by certain scaled B-splines
can be considered as perturbations of the Gabor systems generated by the Gaussian, with a deviation within an arbitrary
small tolerance whenever the order $N$
of the B-spline is sufficiently large. As a consequence we show that
for any choice of translation/modulation parameters $a,b>0$ with $ab<1,$ the scaled version of $B_N$
generates Gabor frames for $N$ sufficiently large.  Considering the
Gabor frame decomposition generated by the Gaussian and a dual window, the results lead to
estimates of the deviation from perfect reconstruction that arise when  the Gaussian  is replaced
by a scaled B-spline, or when the dual window of the Gaussian is replaced by certain explicitly
given and compactly supported linear combinations of the B-splines. In particular, this leads to
a family of approximate dual windows of a very simple form, leading to ``almost perfect reconstruction"
within any desired error tolerance
whenever the product $ab$ is sufficiently small. In contrast, the known (exact) dual windows have a very complicated
form. A similar analysis is sketched
with the scaled B-splines replaced by certain truncations of the Gaussian. As a consequence
of the approach we prove (mostly known) convergence results for the considered scaled B-splines to
the Gaussian in the $L^p$-spaces, as well in the time-domain as in the frequency domain.

\end{abstract}

\begin{minipage}{120mm}


\end{minipage}
\

\section{Introduction} \label{}
For any parameters $a,b\in \mr,$ define the {\it translation operator} $T_a$ and {\it modulation
operators} $E_b$ acting on $\ltr$  by  $T_af(x)= f(x-a), \, E_bf(x)=e^{2\pi i bx} f(x), \, x\in \mr.$
It is well-known (see \cite{Ly,Se2,SW}) that for any $\alpha>0,$ the {\it Gabor system}
$\{E_{mb}T_{na}\varphi\}_{m,n\in \mz}$ generated by the Gaussian $\varphi(x)=e^{-\alpha x^2}$
forms a frame for all $a,b>0$ such that $ab<1.$ The choice of the Gaussian is well motivated
within time-frequency analysis due to the fast decay of the Gaussian and its Fourier transform.
In this paper we consider approximations of the Gaussian via functions with compact support,
namely,  scaled B-splines and truncated Gaussians. We choose to consider
the (centered) B-splines $B_N, \, N\in \mn,$ defined
recursively by
\bes B_1(x)= \chi_{[-1/2, 1/2]}(x), \, \, \, B_{N+1}(x)= B_N * B_1(x)= \int_{-1/2}^{1/2} B_N(x-t)\, dt.\ens

The case of the scaled B-splines is
clearly connected to the central limit theorem; it is also
motivated by results by
Unser et al. \cite{AU}, who proved that certain scalings of the
$N$th order B-splines $B_N$ converge (pointwise
and in certain $L^p$-spaces) to the Gaussian when $N\to \infty.$ Generalizations are
obtained in the papers by Bricks \cite{Bricks} and Goh et al. \cite{GGL}.
We show that for any $a,b>0$ such that
$ab<1,$ the mentioned functions generate Gabor frames  whenever the order of the B-spline
is sufficiently high or the truncation of the Gaussian has sufficiently
large support.  These results are interesting for applications where compact support of the window is desired.

As a consequence of the above results we are able to construct explicitly given approximately dual
frames for Gabor frames $\{E_{mb}T_{na}\varphi\}_{m,n\in \mz}$ generated by the Gaussian, for certain choices of the parameters $a,b>0.$ The approximate dual windows have an extremely simple and explicit form:
they are constant on large intervals, have compact support, and are given explicitly as certain linear
combinations of integer-shifts of the B-splines $B_N$ for some fixed $N\in \mn.$ Furthermore, by choosing the parameter $N$
sufficiently high, we can get as close to perfect reconstruction as desired (see
\eqref{60816a} and the subsequent text for details about this terminology).
  In contrast, the
known (exact) dual windows of the Gaussian  have a very complicated form; see
\cite{J2,Jan9}.

The approach in the current paper leads to proofs of certain other convergence results.
Indeed, we prove that the considered scaled B-splines converge to the Gaussian
in all the $L^p$--spaces for $p\in [1,\infty];$ this result was also
 obtained in \cite{GGL}, while \cite{AU} considered
the case $p\in [2,\infty]$. We also prove that the Fourier transform of the considered scaled
$B_N$ converge uniformly to the Gaussian when $N\to \infty.$

In the rest of this introduction we set the stage by collecting some definitions and necessary results from frame theory.
In Section \ref{060801a} we then prove the mentioned frame results for the scaled B-splines.
In Section \ref{60804a} a similar analysis is sketched for truncated Gaussians.
In Section \ref{60809a} the results are applied to construct explicitly given approximately
dual frames associated with certain Gabor frames generated by the Gaussian.
The consequences
about convergence of the scaled B-splines and its Fourier transform in $L^p$-spaces are collected in Section \ref{60805a}.

A sequence $\{f_k\}_{k\in I}$ in a separable Hilbert space $\h$ is called a {\it frame} if there
exist constants $A,B>0$ such that
\bee \label{70217a} A\, ||f||^2 \le \sum_{k\in I} | \la f, f_k\ra|^2 \le B\,||f||^2, \, \forall f\in \h.\ene
Thus, frames generalize the concept of an orthonormal basis.  The key property of  frames is
that they lead to unconditionally convergent series expansions of the elements in the underlying Hilbert space. Indeed, if
$\{f_k\}_{k\in I}$ is a frame for $\h,$ there always exist a frame $\{g_k\}_{k\in I}$ such that
\bee \label{60816a}  f=\sum_{k\in I} \la f,g_k\ra f_k = \sum_{k\in I} \la f,f_k\ra g_k. \ene
A frame $\{g_k\}_{k\in I}$ satisfying \eqref{60816a} is called a {\it dual frame} of
$\{f_k\}_{k\in I};$ and \eqref{60816a} is expressed by saying that the frames
$\{f_k\}_{k\in I}$ and $\{g_k\}_{k\in I}$ lead to {\it perfect reconstruction.} Note also
that a sequence $\{f_k\}_{k\in I}$ satisfying at least the upper condition in \eqref{70217a}
is called a {\it Bessel sequence with Bessel bound} $B.$

An alternative formulation can be given in terms of operator theory. If the sequence $\{f_k\}_{k\in I}$
satisfies at least the upper frame condition, one defines
the {\it synthesis operator} $T$ by
\bes T: \ell^2(I) \to \h, \, T\{c_k\}_{k\in I} := \sum_{k\in I} c_k f_k.\ens
It is well-known that $T$ is well-defined and bounded.  Denoting the synthesis operator
for the sequence $\{g_k\}_{k\in I}$ by $U,$ perfect reconstruction is equivalent with $TU^*=I.$

A weaker concept was introduced in \cite{CLau}. We say that two Bessel sequences $\{f_k\}_{k\in I}$ and
$\{g_k\}_{k\in I}$ form {\it approximately dual frames} if $||I-TU^*||<1;$ or, in other words,
if there exists a constant $\mu <1$ such that
\bee \label{60822a} \|f-\sum_{k\in I} \la f, {f_k}\ra g_k\|\leq\mu ||f||, \, \forall f\in \h.\ene
The rationale behind approximately dual frames is that all kinds of implementations involve
certain imprecisions; thus, as long as we can control the deviation from perfect reconstruction
measured by the parameter $\mu$ in \eqref{60822a}, approximately dual frames can in practice be
as good as exact dual frames.

Most of the results in the paper will be based on perturbation results, in particular, the result stated below. The result stated in (i) is classic \cite{CH}; we include a short proof of (ii).

\bl \label{pea} Let $\{f_k\}_{k\in I}$ be a frame for a separable Hilbert space $\h$ with bounds $A,B$,
and let $\{\widetilde{f_k}\}_{k\in I}$ be a sequence in $\h$. Then the following hold:
\bei \item[(i)] If there exists a constant
$R<A$ such that \bee \label{2005}
 \sum_{k\in I} |\la f, f_k-\widetilde{f_k}\ra|^2 \le R \ ||f||^2, \ \forall f\in \h,
\ene then $\{\widetilde{f_k}\}_{k\in I}$ is a frame with bounds
$A\left(1- \sqrt{\frac{R}{A}} \right)^2, \
B\left(1+ \sqrt{\frac{R}{B}} \right)^2.$
\item[(ii)]
Assume that  $\{g_k\}_{k\in I}$ is a dual frame of
$\{f_k\}_{k\in I},$ with upper frame bound $B_g$,
and that
\bee \label{60812a}\left\|\sum_{k\in I}c_k(f_k -\widetilde{f_k})\right\|\leq \mu\left(\sum_{k\in I}|c_k|^2\right)^{1/2}\ene
for some $\mu\geq 0$ and all finite sequences $\{c_k\}_{k\in I}$. Then
\bee \label{70220a} \|f-\sum_{k\in I} \la f, \widetilde{f_k}\ra g_k\|\leq\mu \sqrt{B_g}||f||, \, \forall f\in \h.\ene
In particular $\{g_k\}_{k\in I}$ and $\{\widetilde{f_k}\}_{k\in I}$ are approximately dual
frames if $\mu \sqrt{B_g}<1.$
\eni
\el

\noindent {\bf Proof of (ii):} It is well-known that if \eqref{60812a} holds for
all finite sequences, then it holds for all sequences in $\ell^2(I);$ also,
\eqref{60812a} is just a different
way of expressing that $\{f_k-\widetilde{f_k}\}_{k\in I}$ is a
Bessel sequence with bound at most $\mu^2.$
Letting
$c_k= \la f,f_k\ra$ then leads to
\bes \|f-\sum_{k\in I} \la f, \widetilde{f_k}\ra g_k\| & = &
\| \sum_{k\in I} \la f, f_k\ra g_k -\sum_{k\in I} \la f, \widetilde{f_k}\ra g_k\|
\\  & \le & \sqrt{B_g} \left(\sum_{k\in I} |\la f, f_k-\widetilde{f_k}\ra |^2\right)^{1/2} \\
& \leq & \mu \sqrt{B_g}||f||, \, \forall f\in \h,\ens
as desired.\ep
In words, Lemma \ref{pea} (i) says that we can check that a sequence $\{\widetilde{f_k}\}_{k\in I}$ is a frame by proving
that it is ``close" to a frame $\{f_k\}_{k\in I},$ in the sense that $\{f_k-\widetilde{f_k}\}_{k\in I}$ is a
Bessel sequence with a sufficiently small bound.  The importance of Lemma \ref{pea} (ii) lies in the fact that it measures the deviation from
perfect reconstruction which occur when the frame $\{f_k\}_{k\in I}$ is perturbed.

For Gabor systems, the Bessel bound is typically calculated via the so-called CC-condition (see \cite{CC3} or \cite{CB}, Theorem 11.4.2),
stated in Proposition \ref{60804b} (i) below.
Via standard manipulations the CC-condition can
also be formulated in the Fourier domain; since we also need this version of the condition, this
is stated in Proposition \ref{60804b} (ii). We define the Fourier transform of $f\in L^1(\mr)$ by
\bes \widehat{f}(\ga)= \inr f(x) e^{-2\pi i \ga x}dx,\ens with the usual extension to $\ltr.$

\bpr \label{60804b} Let $g\in \ltr,\ a,b>0.$  Then the following hold:

\bei
\item[(i)] If
\bee \label{16g-14}  B(g,a,b) :=\frac1{b} \sup_{x\in [0,a]}
\sum_{k \in \mz} \bigg| \sum_{n\in \mz}  g(x-na)
 \overline{g(x-na- k/b)} \bigg|<\infty, \ene
then $\mts$ is a Bessel sequence with Bessel bound $B(g,a,b)$.
\item[(ii)] If
\bee \widehat{B}(g,a,b) :=\frac1{a}
\sup_{\ga \in [0,b]} \sum_{k \in \mz}  \bigg| \sum_{n\in \mz}
 \widehat g(\ga-nb)\overline{\widehat g(\ga-nb- k/a)} \bigg|
\label{16g-0}<\infty,\ene
then $\mts$ is a Bessel sequence with Bessel bound $\widehat{B}(g,a,b)$.
\eni
\epr

\section{Approximation of the Gaussian via B-splines} \label{060801a}
Our  first goal is to show
that for any choice of the translation/modulation parameters, the Bessel bound of the Gabor system generated by
the difference between certain scaled B-splines and the Gaussian tend to zero as $N\to \infty.$
Eventually this leads to the main conclusion in this section, namely, that for any $a,b>0$ with
$ab<1,$ the scaled B-spline $B_N$ generate Gabor frames whenever the order $N$ is sufficiently large.

In order to estimate the difference between a Gaussian and the scaled B-splines
pointwise, we will need the following result concerning the pointwise difference between
a Gaussian and the $N$th power of the sinc-function. We will state two versions: a precise
pointwise estimate, and a qualitative estimate, concerning the behavior as $N\to \infty.$
Given a function $f: \mn \to \mc$ the notation $f(N)=o(1)$ means
that for every $\epsilon>0$, there exists $N_0\in\mn$ such that
for $N\geq N_0$, $|f(N)|\leq \epsilon.$ Even when several parameters are involved
like in \eqref{70217d} below, the symbol $o(1)$ will always refer to the dependence on the parameter $N.$

\bl\label{16g-1} Let $N\geq 14$ be an integer. Then
\begin{equation} \label{70217d}
\left| e^{-\frac{x^2 N}{6}}-\left(\frac{\sin x}{x}
\right)^N \right|
\leq     \left\{
    \begin{array}{ll}
\frac{4}{5e^{2}N}\left(1+\frac{17\ln N}{7N}\right)=\frac{4}{5e^{2}N}(1+o(1)), & |x|< \sqrt{\frac{12 \ln N}{N}}; \\
\frac{4(\ln
N)^2}{5N^3}\left(1 +\frac{17\ln N}{7N}\right)=\frac{4(\ln N)^2}{5N^3}\left(1 +o(1)\right), & \sqrt{\frac{12 \ln N}{N}}\leq |x|\leq \frac{\pi}{2}.
    \end{array}
    \right.
\end{equation}

\el
\bp Fix $N\ge 14.$ By symmetry it is enough to consider $x\in [0,\frac{\pi}{2}]$.  A crucial step in the proof is to observe that
\bee \label{60804c} e^{-\frac{x^2}{6}}\geq \frac{\sin x}{x}, \ \forall x\in [0,\frac{\pi}{2}],\ene so let us argue for this first.
Taylor's theorem implies that
$$\frac{\sin x}{x} \leq 1-\frac{x^2}{3!}+\frac{x^4}{5!}  \hspace{.3cm}
\mbox{and}  \hspace{.3cm}
e^{-\frac{x^2}{6}} \geq 1-\left(\frac{x^2}{6}\right)+\frac{1}{2!}\left(\frac{x^2}{6}\right)^2
-\frac{1}{3!}\left(\frac{x^2}{6}\right)^3. $$
A direct calculation shows that
$1-\left(\frac{x^2}{6}\right)+\frac{1}{2!}\left(\frac{x^2}{6}\right)^2
-\frac{1}{3!}\left(\frac{x^2}{6}\right)^3 \geq 1-\frac{x^2}{3!}+\frac{x^4}{5!}.$
Thus we obtain that \eqref{60804c} holds for $x\in [0,\frac{\pi}{2}],$  as desired.

Now we will use that
\begin{equation}\label{16g-11}
\left|e^{-\frac{x^2 N}{6}}-\left(\frac{\sin x}{x} \right)^N
\right|= e^{-\frac{x^2 N}{6}} \left|1- \left( \frac{\sin x}{x}\,
e^{\frac{x^2}{6}} \right)^N \right|.
\end{equation}
Fix $x\in [0,\frac{\pi}{2}]$. Taylor's theorem
again implies that $\frac{\sin x}{x} \geq 1-\frac{x^2}{3!}+\frac{x^4}{5!}-\frac{x^6}{7!}$ and
$e^{\frac{x^2}{6}} \geq 1+\frac{x^2}{6}+\frac{x^4}{72}$. In the following
calculation we will use that
$\frac{x^8}{12096}-\frac{x^{10}}{362880}=\frac{x^8}{362880}(30-x^2)>0$.
Indeed, this implies that
\begin{eqnarray*}
\left(\frac{\sin x}{x}\, e^{\frac{x^2}{6}}\right)^N
&\geq&
\left(1-\frac{x^2}{3!}+\frac{x^4}{5!}-\frac{x^6}{7!} \right)^N
\left(1+\frac{x^2}{6}+\frac{x^4}{72}\right)^N \\
&=&\left(1-\frac{x^4}{180}-\frac{17x^6}{15120}+\frac{x^8}{12096}-\frac{x^{10}}{362880}\right)^N\\
&\geq&\left(1-\frac{x^4}{180}-\frac{17x^6}{15120}\right)^N.
\end{eqnarray*}

Using the inequality $(1-t)^N\geq
1-Nt$ for $0<t<1$  and that  $\frac{x^4}{180}+\frac{17x^6}{15120}<1$, we have
$\left(\frac{\sin x}{x}\, e^{\frac{x^2}{6}}\right)^N\geq
1-N(\frac{x^4}{180}+ \frac{17x^6}{15120})=1-\frac{N x^4}{180}(1+\frac{17}{84}x^2).$
This together with \eqref{60804c} implies that
$0\leq 1-\left(\frac{\sin x}{x}\, e^{\frac{x^2}{6}}\right)^N \leq
\frac{N x^4}{180}(1+\frac{17}{84}x^2). $
Therefore  we have
\begin{equation}\label{16g-38}
\left|e^{-\frac{x^2 N}{6}}-\left(\frac{\sin x}{x} \right)^N
\right|\leq e^{-\frac{x^2 N}{6}}\, \frac{N x^4}{180}(1+\frac{17}{84}x^2)=:f_N(x),
\end{equation}
by \eqref{16g-11}.
Note that $\sqrt{\frac{12\ln N}{N}}\leq \frac{\pi}{2}$ for $N\geq
14$.
We now split $[0,\frac{\pi}{2}]$ into two intervals
$[0,\sqrt{\frac{12\ln N}{N}}[$ and $[\sqrt{\frac{12\ln
N}{N}},\frac{\pi}{2}]$:

(1) For $x\in [0,\sqrt{\frac{12\ln N}{N}}[$,
\eqref{16g-38} implies that
$$\left|e^{-\frac{x^2 N}{6}}-\left(\frac{\sin x}{x} \right)^N \right|\leq
e^{-\frac{x^2 N}{6}}\, \frac{N x^4}{180}\left(1+\frac{17\ln N}{7N}\right).$$
Fixing $N\in\mn$, the function $x\mapsto e^{-\frac{x^2 N}{6}}\, \frac{N x^4}{180}$
attains its maximum value $\frac{4}{5e^{2}N}$ on the interval $[0,\sqrt{\frac{12\ln N}{N}}[$
at $x=\sqrt{\frac{12}{N}}$; thus
\begin{equation*}
\left|e^{-\frac{x^2 N}{6}}-\left(\frac{\sin x}{x} \right)^N \right|
\leq \frac{4}{5e^{2}N}\left(1+\frac{17\ln N}{7N}\right) =
\frac{4}{5e^{2}N}(1+o(1)).
\end{equation*}

(2) Note that for $N\ge 14$, the function $f_N$ is decreasing
on the interval $[\sqrt{\frac{12\ln N}{N}}, \frac{\pi}{2}]$.
Thus, for $x\in [\sqrt{\frac{12\ln N}{N}}, \frac{\pi}{2}]$,
\begin{equation*}
f_N(x) \leq  f_N\left(\sqrt{\frac{12\ln N}{N}}\right)
= \frac{4(\ln
N)^2}{5N^3}\left(1 +\frac{17\ln N}{7N}\right)
=\frac{4(\ln
N)^2}{5N^3}\left(1 +o(1)\right),
\end{equation*}
as desired. \ep

In the next lemma we introduce the exact scaling of the B-splines that we use throughout
the paper; it goes back to the paper \cite{AU} by Unser et al. The result describes the asymptotic
behavior of these functions.



\bl \label{60802c}
Let $N\in\mn$ and let
\bee \label{60802d} p_N(x):=\frac{1}{\sqrt{2\pi}}e^{-x^2/2} -
\sqrt{\frac{N}{12}} \, B_N\left(\sqrt{\frac{N}{12}}\,
x\right).\ene Given any $a,b>0,$ the Bessel bound $ \widehat{B}p_N,a,b)$ in
\eqref{16g-0} can be estimated by

\begin{equation}\label{16g-22a}
\widehat{B}(p_N,a,b) \leq
 \frac{64\ln N}{b25 \pi^2 e^4 N^2}
\left(1+o(1)\right).
\end{equation}
In particular, for  any $\epsilon>0$, there exists a positive integer $N(a,b)$ such that
for $N\ge N(a,b)$ the Gabor system $\{E_{mb}T_{na}p_N\}_{m,n\in \mz}$ is
a  Bessel sequence with bound at most $\epsilon.$
\el

\bp
Let $N\geq 14$ be an integer. Using that $\widehat{B_N}(\ga)= \left(\frac{\sin \pi \ga}{\pi \ga}\right)^N,$ a change of variable
yields that
\begin{equation}\label{16g-12}
\widehat{p_{N}}(\ga) = e^{-2\pi^2 \ga^2}
-\left(\frac{\sin\left(\pi\sqrt{\frac{12}{N}}\, \ga \right)}{\pi\sqrt{\frac{12}{N}}\,\ga}\right)^{N}.
\end{equation}
By letting $x=\pi \sqrt{\frac{12}{N}}\, \ga$,
Lemma \ref{16g-1} implies that
\begin{equation}\label{16g-7}
\left| e^{-2\pi^2 \ga^2}
-\left(\frac{\sin\left(\pi\sqrt{\frac{12}{N}}\, \ga \right)}{\pi\sqrt{\frac{12}{N}}\,\ga}\right)^{N}\right|
\leq
\left\{
\begin{array}{ll}
\frac{4}{5e^{2}N}\left(1+\frac{17\ln N}{7N}\right),&  |\ga|< \frac{\sqrt{\ln N}}{\pi};\\
\frac{4(\ln N)^2}{5N^3}\left(1+\frac{17\ln N}{7N}\right),&  \frac{\sqrt{\ln
N}}{\pi}\leq |\ga|\leq \frac{\sqrt{N}}{4\sqrt{3}}.
\end{array}
\right.
\end{equation}
Also, clearly
\begin{equation}\label{16g-8}
\left| e^{-2\pi^2 \ga^2}
-\left(\frac{\sin\left(\pi\sqrt{\frac{12}{N}}\, \ga \right)}{\pi\sqrt{\frac{12}{N}}\,\ga}\right)^{N}\right|
\leq e^{-2\pi^2 \ga^2} +
\left|\frac{1}{\pi\sqrt{\frac{12}{N}}\,\ga}\right|^N, \ \ga\in\mr.
\end{equation}
We now estimate the Bessel bound \eqref{16g-0} for $g=p_N.$  We first note that
\bee \notag \widehat{B}(p_N, a,b)  & \leq &\frac1{a}
\sup_{\ga \in [0,b]} \sum_{k \in \mz}  \sum_{n\in \mz}
\bigg| \widehat{p_N}(\ga-nb)\widehat{p_N}(\ga-nb- k/a) \bigg| \\ & \le &
\frac1{a}\left(\sup_{\ga \in [0,b], n\in \mz} \sum_{k \in \mz}
\bigg|\widehat{p_N}(\ga-nb- k/a) \bigg|\right) \left(\sup_{\ga \in [0,b]} \sum_{n\in \mz}
\bigg| \widehat{p_N}(\ga-nb)\bigg|\right). \label{60802a} \ene
\noindent {\bf Estimate of
$\sum_{k \in \mz}
\bigg| \widehat{p_{N}}(\ga-nb- k/a) \bigg|$ :}
 For a fixed $n\in\mz$, $\ga-nb- k/a$ hits the interval
$]-\frac{\sqrt{\ln N}}{\pi},\frac{\sqrt{\ln N}}{\pi}[$ for at most
$\lceil \frac{2a\sqrt{\ln N}}{\pi} \rceil$
 values of $k\in \mz$ and its
contribution is at most
$\lceil \frac{2a\sqrt{\ln N}}{\pi} \rceil
\frac{4}{5e^{2}N}(1+\frac{17\ln N}{7N})$
 by \eqref{16g-7}. Also,
by \eqref{16g-7}, the contribution from $\ga-nb- k/a$ hitting the
interval $]-\frac{\sqrt{N}}{4\sqrt{3}},-\frac{\sqrt{\ln
N}}{\pi}]\cup [\frac{\sqrt{\ln
N}}{\pi},\frac{\sqrt{N}}{4\sqrt{3}}[$ is at most
$ 2 \lceil a ( \frac{\sqrt{N}}{4\sqrt{3}}-\frac{\sqrt{\ln
N}}{\pi} )\rceil \frac{4(\ln N)^2}{5N^3}(1+\frac{17\ln
N}{7N})$. By \eqref{16g-8},
the contribution from $\ga-nb- k/a$ hitting the interval
$[\frac{\sqrt{N}}{4\sqrt{3}},\infty[$ is at most
\begin{eqnarray}
\sum_{k=0}^{\infty} \left| \widehat{p_{N}}
\left(\frac{\sqrt{N}}{4\sqrt{3}}+\frac{k}{a}\right)\right| &\leq &\sum_{k=0}^{\infty } e^{-2\pi^2
(\frac{\sqrt{N}}{4\sqrt{3}}+\frac{k}{a})^2} +
\sum_{k=0}^{\infty } \left|\frac{\sqrt{N}}
{\sqrt{12} \pi \left( \frac{\sqrt{N}}{4\sqrt{3}}+\frac{k}{a} \right)} \right|^{N} \nonumber\\
&=:&(I)+(II)
\end{eqnarray}
We now estimate (I) and (II) as follows:
\begin{eqnarray*}
 (I)\leq e^{-\frac{\pi^2 N}{24}}  \sum_{k=0}^{\infty} e^{-\left(\frac{\pi^2 \sqrt{N}}{\sqrt{3}a}\right)k}
= \frac{e^{-\frac{\pi^2 N}{24}}}{1-e^{-\frac{\pi^2 \sqrt{N}}{\sqrt{3}a}} }
\end{eqnarray*}
and
\begin{eqnarray*}
   (II)&=&
\sum_{k=0}^{\infty } \left|
\frac{\sqrt{N}}{\sqrt{12} \pi \frac{\sqrt{N}}{4\sqrt{3}}
\left(1 +\frac{4\sqrt{3}k}{a\sqrt{N}} \right)}
\right|^{N}   \\
&=&\left( \frac{2}{\pi} \right)^{N}\sum_{k=0}^\infty
\left(1+\frac{ 4\sqrt{3}k}{a\sqrt{N}}\right)^{-N}\\
   &\leq& \left( \frac{2}{\pi} \right)^{N}
   \left(1 +\int_0^\infty \left(1+\frac{4 \sqrt{3}x}{a\sqrt{N}} \right)^{-N} dx \right) \\
   &=& \left( \frac{2}{\pi} \right)^{N} \left(1+\frac{a\sqrt{N}}{4\sqrt{3}(N-1)} \right).
\end{eqnarray*}
The contribution from $]-\infty,-\frac{\sqrt{N}}{4\sqrt{3}}]$ is at most the same as above.
Therefore,



\begin{eqnarray}
\sum_{k \in \mz}
\bigg| \widehat{p_{N}}(\ga-nb- k/a) \bigg|
&\leq&  \left(1+\frac{2 a\sqrt{\ln N}}{\pi} \right)\frac{4}{5e^{2}N}\left(1+\frac{17\ln N}{7N}\right) \nonumber \\
&& +\left(2+2 a\left( \frac{\sqrt{N}}{4\sqrt{3}}-\frac{\sqrt{\ln N}}{\pi} \right)\right)
\frac{4(\ln N)^2}{5N^3}\left(1+\frac{17\ln N}{7N}\right)  \nonumber\\
&&+
2 \left(\frac{e^{-\frac{\pi^2 N}{24}} }{1-e^{-\frac{\pi^2 \sqrt{N}}{\sqrt{3}a}} }
+\left( \frac{2}{\pi} \right)^{N} \left(1+\frac{a\sqrt{N}}{4\sqrt{3}(N-1)} \right)\right) \nonumber\\
&=:& P_N(a) \label{16g-21} \\
&=& \frac{8 a\sqrt{\ln N}}{5\pi e^{2}N} (1+o(1)). \nonumber
\end{eqnarray}
Note that the above estimate  is independent on the choice of $n\in \mz$ and $b>0.$

\noindent {\bf Estimate of
$\sum_{n \in \mz}
\bigg| \widehat{p_{N}}(\ga-nb) \bigg|$ :}  Using that  the estimate \eqref{16g-21}
also holds in the particular case where the parameter $a$ takes the value $1/b,$
it follows  that
\begin{eqnarray*}
\sum_{n \in \mz}
\bigg| \widehat{p_{N}}(\ga-nb) \bigg|
&\leq&
\left(1+ \frac{2\sqrt{\ln N}}{b\pi}  \right)\frac{4}{5e^{2}N}\left(1+\frac{17\ln N}{7N}\right) \\
&&+\left(2+ \frac{2}{b}\left( \frac{\sqrt{N}}{4\sqrt{3}}-\frac{\sqrt{\ln N}}{\pi} \right) \right)
\frac{4(\ln N)^2}{5N^3}\left(1+\frac{17\ln N}{7N}\right)
 \\
&& + 2 \left(
\frac{e^{-\frac{\pi^2 N}{24}} }{1-e^{-\frac{\pi^2 b\sqrt{N}}{\sqrt{3}}} }
+\left( \frac{2}{\pi} \right)^{N} \left(1+\frac{\sqrt{N}}{4b\sqrt{3}(N-1)} \right)\right)\\
&=&
P_N(1/b)=
 \frac{8\sqrt{\ln N}}{b5\pi e^{2}N} (1+o(1)).
\end{eqnarray*}
Hence \eqref{60802a} implies that the Bessel bound of $\{E_{mb}T_{na}p_N\}_{m,n\in \mz}$
can be estimated by
\begin{equation}\label{16g-22}
\widehat{B}(p_N,a,b) \leq \frac{ P_N(a) P_N(1/b)}{a}= \frac{64\ln N}{b25 \pi^2 e^4 N^2}
\left(1+o(1)\right).
\end{equation}
It is now clear that
$\widehat{B}(p_N,a,b) \rightarrow 0$ as $N\rightarrow \infty$.
\ep

Lemma \ref{60802c} now allows  us to give an easy proof for the
main result in this section, namely, that for any $a,b>0$ with $ab<1,$ the considered scaled B-splines $B_N$ generate Gabor frames
whenever the order $N$ is sufficiently high.

\bt \label{16g-17}
For any $a,b>0$ such that $ab<1,$ there exists a positive integer $N(a,b)$
such that the function
\bee \label{60802g} g_N(x):=  \sqrt{\frac{N}{12}} \, B_N\left(\sqrt{\frac{N}{12}}\, x\right)\ene
generates a Gabor frame $\{E_{mb}T_{na}g_N\}_{m,n\in\mz}$   whenever $N\ge N(a,b).$ \et
\bp We have already mentioned that the Gaussian generates a Gabor frame
for  any $a,b>0$ such that $ab<1.$  Fixing such $a,b,$ let $A$ denote a lower
frame bound for the Gabor frame generated by the function $\frac{1}{\sqrt{2\pi}}e^{-x^2/2}.$
By Lemma \ref{60802c}
we can choose an integer $N(a,b)$ such that the Bessel bound $\widehat{B}(p_N, a,b)$ for
the function $p_N$ in \eqref{60802d} is smaller than $A$ for $N\ge N(a,b).$ By Lemma
\ref{pea} this implies that the function $g_N$ in \eqref{60802g} generates a
Gabor frame for $N\ge N(a,b).$ \ep

Note that the frame properties of the (non-scaled) B-splines  is a very active research area.
We will not directly use any of the results from the literature, so we just refer to the
recent paper \cite{LemNiel} and the extensive list of references therein for more information.

The proof of Theorem \ref{16g-17} gives some information about how to choose the integer $N(a,b):$
it should be chosen such that the Bessel bound $\widehat{B}(p_N, a,b)$ for
the function $p_N$ in \eqref{60802d} is smaller than the lower frame bound $A$ for the Gabor frame generated
by the Gaussian.  The next result, which is a consequence
 of the calculations in the proof of Lemma \ref{60802c}, gives more explicit information about how to choose $N(a,b)$
such that the Bessel bound for the function $p_N$ does not exceed a prescribed maximal value:

\bc
Let $N_0\geq 14$ be an integer, and let
\begin{eqnarray}
K(x)&=&  \left(  \frac{4}{5e^2} \left(\frac{2 x }{\pi}+\frac{1}{\sqrt{\ln N_0}}\right)
+\left(\frac{\ln N_0}{N_0}\right)^{3/2} \frac{4}{5 }
\left( \frac{2}{\sqrt{N_0}}+\frac{x}{2\sqrt{3}} \right)\right) \left(1+\frac{17\ln N_0}{7N_0}\right) \nonumber\\
&&+ \frac{2N_0}{\sqrt{\ln N_0}}\left(\frac{e^{-\frac{\pi^2 N_0}{24}} }{1-e^{-\frac{\pi^2 \sqrt{N_0}}{\sqrt{3}x}} }
+\left( \frac{2}{\pi} \right)^{N_0} \left(1+\frac{x\sqrt{N_0}}{4\sqrt{3}(N_0-1)} \right)\right). \label{16g-35}
\end{eqnarray}
Define $p_{N}$ as in \eqref{60802d}.
Given any $\epsilon>0$, choose
\begin{equation}\label{16g-27}
N(a,b):=\left\lfloor \frac{K(a)K(1/b)}{a\epsilon} \right\rfloor +N_0.
\end{equation}
Then     $\widehat{B}(p_{N(a,b)},a,b)< \epsilon.$
\ec
\bp
First, fix any integer $N_0\geq 14.$
It is easy to show by induction that
for $\alpha\leq  \frac{14}{15}$ and $N\geq N_0$, we have that $\alpha^N \leq \frac{\alpha^{N_0}N_0}{N}.$
Since $e^{-\pi^2/24}\leq \frac{14}{15}$ and $2/\pi \leq \frac{14}{15}$
 we can now use this
to conclude that
\begin{equation}
e^{-\frac{\pi^2 N}{24}} \leq
 e^{-\frac{\pi^2 N_0}{24}}\frac{N_0}{N} \text{ and }
\left( \frac{2}{\pi} \right)^{N} \leq
\left( \frac{2}{\pi} \right)^{N_0} \frac{N_0 }{N }.
\end{equation}
Moreover, the functions $N\mapsto \frac{\ln N}{N}$,
$N\mapsto e^{-\frac{\pi^2 \sqrt{N}}{a\sqrt{3}}}$,
$N\mapsto \frac{\sqrt{N}}{N-1}$  are decreasing.
Thus, with the given definition of $K(x)$ and $P_N$ defined by \eqref{16g-21},
\begin{eqnarray}
P_N(a)
&\leq&
\frac{\sqrt{\ln N}}{N}
\left(  \frac{4}{5e^2} \left(\frac{2 a }{\pi}+\frac{1}{\sqrt{\ln N}}\right)
\left(1+\frac{17\ln N}{7N}\right)\right. \nonumber\\
&&\left.
+ \left(\frac{\ln N}{N}\right)^{3/2} \frac{4}{5 } \left( \frac{2}{\sqrt{N}}+\frac{a}{2\sqrt{3}} \right)
 \left(1+\frac{17\ln N}{7N}\right) \right. \nonumber\\
&&\left.
+2 \frac{N}{\sqrt{\ln N}}\left(\frac{e^{-\frac{\pi^2 N}{24}} }{1-e^{-\frac{\pi^2 \sqrt{N}}{\sqrt{3}a}} }
+\left( \frac{2}{\pi} \right)^{N} \left(1+\frac{a\sqrt{N}}{4\sqrt{3}(N-1)} \right)\right)\right) \nonumber\\
&\leq&\frac{\sqrt{\ln N}}{N}K(a).   \label{16g-26}
\end{eqnarray}
This together with \eqref{16g-22} implies that
$$\widehat{B}(p_N,a,b)\leq  \frac{K(a)K(1/b)\ln N}{a N^2}\leq
 \frac{K(a)K(1/b)}{aN}.$$
 Note that we have that
$\frac{K(a)K(1/b)}{a\epsilon} < N(a,b)$ by \eqref{16g-27}.
Thus
$\widehat{B}(p_{N(a,b)},a,b)\leq \frac{K(a)K(1/b)}{aN(a,b)} < \epsilon,$ as desired. \ep

\section{Frame properties of truncated Gaussians} \label{60804a}

In this section we will repeat the procedure of Section \ref{060801a} with the scaled
B-splines replaced by certain compactly supported  ``truncations" of the Gaussian,
namely, functions of the form
\bes \varphi_N(x):= \left(g(x)-g(N)\right) \chi_{[-N,N]}(x)\ens
for some $N\in \mn$ and with $g(x)=e^{-x^2/2}.$ Note that the subtraction of the value $g(N)$
implies that $\varphi_N$ is continuous.
Due to the fast decay of the exponential function $g$, the function $\varphi_N$ will in general provide a very
good approximation of $g$ for large values of $N.$
In the next lemma we show that the truncated versions of the Gaussian approximate the Gaussian well
in terms of the Bessel condition, whenever $N$ is large. Similar considerations can
of course be done for other functions, but for the sake of the flow of the current paper
we stick to the Gaussian.

\bl \label{60809b} Let $N\in \mn$.
Consider
\begin{equation}\label{16g-28}
q_N(x):= e^{-x^2/2}-\left(e^{-x^2/2} -e^{-N^2/2} \right)\chi_{[-N,N]}(x).
\end{equation}
Given any $a,b>0,$ the Bessel bound $ B(q_N,a,b)$ in \eqref{16g-14} can be estimated by
\begin{eqnarray*} 
B(q_N,a,b) &\leq&
\left(1+2bN+\frac{2}{1- e^{-\frac{N}{b}}} \right)
\left( 1+\frac{2N}{a}+\frac{2}{1-e^{-aN}} \right)\frac{e^{-N^2}}{b}\\
&=& \frac{4N^2e^{-N^2}}{a}\left(1+o(1)\right).
\end{eqnarray*}
In particular, for  any $\epsilon>0$, there exists a positive integer $N(a,b)$ such that
for $N\ge N(a,b)$ the Gabor system $\{E_{mb}T_{na}q_N\}_{m,n\in \mz}$ is
a  Bessel sequence with bound at most $\epsilon.$
\el

\bp The proof is similar to the proof of Lemma \ref{60802c}, so we only sketch it.
We estimate the Bessel bound \eqref{16g-14} for $g=q_N.$  We note that
\begin{equation}\label{60802a-1}
 B(q_N, a,b)  \le
\frac1{b}\left(\sup_{x \in [0,a], n\in \mz} \sum_{k \in \mz}
\bigg|q_N(x-na- k/b) \bigg|\right)
\left(\sup_{x \in [0,a]} \sum_{n\in \mz}
\bigg| q_N(x-na)\bigg|\right).
\end{equation}
We also note that
\begin{equation}\label{16g-15}
  q_N(x)=\left\{
  \begin{array}{ll}
  e^{-N^2/2}, & x\in[-N,N];\\
  e^{-x^2/2}, & x\in[-N,N]^c.
  \end{array}
  \right.
\end{equation}
\noindent {\bf Estimate of
$\sum_{k \in \mz}
\bigg|  q_{N}(x-na- k/b) \bigg|$ :}
For a fixed $n\in\mz$, $x-na- k/b$ hits
$]-N,N[$
at most
$\lceil 2bN \rceil$
 times for $k$'s and
its contribution is at most
$\lceil 2bN \rceil e^{-N^2/2}$
 by \eqref{16g-15}.
By \eqref{16g-15} again, the contribution from $x-na- k/b$ hitting the interval
$[N,\infty[$ is at most
\begin{eqnarray*}
\sum_{k=0}^{\infty} \left|  q_N
\left(N+\frac{k}{b}\right)\right|
\leq \sum_{k=0}^{\infty }
e^{-\frac{N^2}{2} -\frac{Nk}{b}} = \frac{e^{-\frac{N^2}{2}}}{1- e^{-\frac{N}{b}}}.
\end{eqnarray*}
The contribution from $]-\infty,-N]$ is at most the same as above.
Therefore,
\begin{equation}\label{16g-18}
\sum_{k \in \mz}
\bigg|  q_{N}(x-na- k/b) \bigg|
\leq(1+ 2bN) e^{-\frac{N^2}{2}}+ \frac{2 e^{-\frac{N^2}{2}}}{1- e^{-\frac{N}{b}}}
=:Q_N(b)= 2bN e^{-\frac{N^2}{2}}(1+o(1)).
\end{equation}
\noindent {\bf Estimate of
$\sum_{n \in \mz}
\bigg|  q_{N}(x-na) \bigg|$ :} Applying  the above result
with $b=1/a$
yields  that
\begin{equation}\label{16g-19}
\sum_{n \in \mz}
\bigg|  q_{N}(x-na) \bigg|
\leq
(1+\frac{2N}{a}) e^{-\frac{N^2}{2}}+ \frac{2e^{-\frac{N^2}{2}}}{1- e^{-aN}}
=Q_N(1/a)= \frac{2Ne^{-\frac{N^2}{2}} }{a}(1+o(1)).
\end{equation}
Hence \eqref{60802a-1} together with \eqref{16g-18} and \eqref{16g-19}
implies that the Bessel bound of $\{E_{mb}T_{na}q_N\}_{m,n\in \mz}$
can be estimated by
\begin{eqnarray}
B(q_N,a,b) &\leq& \frac{Q_N(b)Q_N(1/a)}{b} \nonumber \\
&=& \left(1+2bN+\frac{2}{1- e^{-\frac{N}{b}}} \right)
\left( 1+\frac{2N}{a}+\frac{2}{1-e^{-aN}} \right)\frac{e^{-N^2}}{b}\nonumber\\
&=&\frac{4N^2e^{-N^2}}{a}\left(1+o(1)\right). \label{16g-16}
\end{eqnarray}
This clearly implies that
$B(q_N,a,b) \rightarrow 0$ as $N\rightarrow \infty$, as desired.
\ep

As in  Theorem \ref{16g-17}, Lemma \ref{60809b} has the following consequence:
\bt \label{060824a} For any $a,b>0$ such that $ab<1,$ there exists a positive integer $N(a,b)$
such that the function
\bee \label{60805b}
\varphi_N(x):=\left(e^{-x^2/2} -e^{-N^2/2} \right)\chi_{[-N,N]}(x)
\ene
generates a Gabor frame $\{E_{mb}T_{na} \varphi_N\}_{m,n\in\mz}$ whenever $N\ge N(a,b).$ \et

Analog to the case of B-spline approximations in Section \ref{060801a} one can
derive a quantitative estimate on the number $N(a,b)$ in Theorem \ref{060824a}; we leave the details
to the interested reader.

\section{Approximately dual frames} \label{60809a}
The known dual windows associated with the Gabor frames generated by the Gaussian are
quite complicated  (see the end of the section for a more detailed discussion).
On the other hand, for sufficiently small
modulation parameters $b>0$ and the translation parameter $a=1,$ B-spline generated Gabor frames
$\{E_{mb}T_nB_N\}{m,n\in\mz}$
have convenient dual windows that are finite linear combinations of integer-translates
of the B-spline \cite{C19,CR}. In this section we prove that replacing the exact dual window of the
Gaussian by these simple dual windows of the B-spline can get us as close to perfect reconstruction as
desired.

The starting point is the observation  that  for certain choices of translation/modulation parameters one can construct
a dual frame associated with the Gabor frames generated by the considered scaled B-splines; see
\cite{C19, CR} for details.

\begin{figure}
\centerline{
\includegraphics[width=5in,height=1.3in]{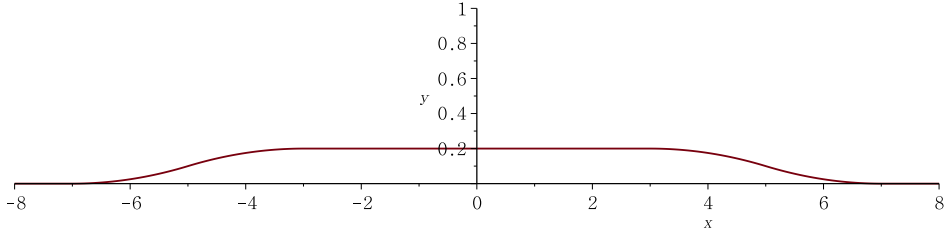}}\hfil
\caption{The dual window $h_N$ in Proposition \ref{16g-33} for $N=3$.                    } \label{16g-36}
\end{figure}

\bpr \label{16g-33}
Let $N\in\mn$. Let $a:=\sqrt{\frac{12}{N}}$ and $b:=\frac{1}{2N-1}\sqrt{\frac{N}{12}}$.
Then the functions
\begin{equation}\label{16g-30}
g_N(x):=\sqrt{\frac{N}{12}} B_N\left(\sqrt{\frac{N}{12}}\, x\right)
\end{equation}
and
\begin{equation}\label{16g-31}
h_N(x):=b \sqrt{\frac{12}{N}} \sum_{n=-N+1}^{N-1} B_N\left(\sqrt{\frac{N}{12}}\, x+n\right)
\end{equation}
generate a pair of dual frames $\{E_{mb}T_{na}\, g_N\}_{m,n\in\mz}$ and
$\{E_{mb}T_{na}\, h_N\}_{m,n\in\mz}$.
\epr
Note that in Proposition \ref{16g-33} we have $ab= (2N-1)^{-1};$ thus, for large values of
$N\in \mn$ the constructed Gabor frames are highly redundant.
Figure \ref{16g-36} shows the dual window $h_N$  for $N=3$. The shape of $h_3$ is typical for
all the functions $h_N, \, N\in \mn;$ indeed, $h_N$ is constant on the interval $[-\sqrt{3N},\sqrt{3N}]$ and
supported on the interval $[-(\frac{3N}{2}-1)\sqrt{\frac{12}{N}},(\frac{3N}{2}-1)\sqrt{\frac{12}{N}}].$

\bl
Let $N\in\mn$. Let $a:=\sqrt{\frac{12}{N}}$
and $b:=\frac{1}{2N-1}\sqrt{\frac{N}{12}}$.
Define $h_N$ as in  \eqref{16g-31}.
Then
the Bessel bound $B(h_N,a,b)$ in \eqref{16g-14} for $\{E_{mb}T_{na}\, h_N\}_{m,n\in\mz}$
can be estimated by
\begin{equation} \label{16g-32}
 B\left(h_N, a, b\right)
\leq
\frac{9\sqrt{3}}{2\sqrt{N}}(1+o(1)).
\end{equation}
\el
\bp We estimate the Bessel bound \eqref{16g-14} for $g=h_N$
with $a=\sqrt{\frac{12}{N}}$ and
$b=\frac{1}{2N-1}\sqrt{\frac{N}{12}}$.
As we already saw in  \eqref{60802a-1},
\begin{equation}\label{60802a-2}
 B(h_N, a,b)  \le
\frac1{b}\left(\sup_{x \in [0,a], n\in \mz} \sum_{k \in \mz}
\bigg|h_N(x-na- k/b) \bigg|\right)
\left(\sup_{x \in [0,a]} \sum_{n\in \mz}
\bigg| h_N(x-na)\bigg|\right).
\end{equation}
Since $B_N$ satisfies the partition of unity property, it follows directly from
our knowledge of the support of $h_N$ that
$$h_N(x)\leq b\sqrt{\frac{12}{N}}
\chi_{\left[-\left(\frac{3N}{2}-1\right)\sqrt{\frac{12}{N}},\left(\frac{3N}{2}-1\right)\sqrt{\frac{12}{N}}\right]}(x).$$
By arguments similar to the ones used in the proof of Lemma \ref{60809b}, we have
$$ \sum_{k \in \mz}
\bigg|  h_{N}(x-na- k/b) \bigg|\leq
\left(1+2b\left(\frac{3N}{2}-1\right)\sqrt{\frac{12}{N}}\right)b\sqrt{\frac{12}{N}}
$$
and
$$\sum_{n \in \mz}
\bigg|  h_{N}\left(x-na\right) \bigg|\leq
\left(1+\frac{2}{a}\left(\frac{3N}{2}-1\right)\sqrt{\frac{12}{N}}\right)b\sqrt{\frac{12}{N}}
.$$
By \eqref{60802a-2} and direct calculation based on the above estimates
and the given values of $a,b$, the Bessel bound of
$\{E_{mb}T_{na}\, h_N\}_{m,n\in \mz}$
can be estimated by
\begin{eqnarray*}
  B\left(h_N, \sqrt{\frac{12}{N}}, \frac{1}{2N-1}\sqrt{\frac{N}{12}}\right)
\leq
\frac{9\sqrt{3}}{2\sqrt{N}}(1+o(1)),
\end{eqnarray*}
as desired.
\ep

The next result proves that for $N\in \mn$ sufficiently large and certain choices of $a,b>0,$
the Gabor system $\{E_{mb}T_{na}h_N\}_{m,n\in \mz}$ is an approximately dual frame of
a Gabor system generated by a (multiple of a) Gaussian:

\bt \label{60812f}
 Let $N\in \mn$, and let $a:=\sqrt{\frac{12}{N}}$ and $b:=\frac{1}{2N-1}\sqrt{\frac{N}{12}}$.
Define $h_N$ as in \eqref{16g-31} and let
$$\widetilde g(x)=\frac{1}{\sqrt{2\pi}}\, e^{-x^2/2}.$$
Then
\begin{equation}\label{16g-34}
\left\| f-\sum_{m,n\in\mz}\langle f, E_{mb}T_{na} \widetilde g \rangle
E_{mb}T_{na} h_N \right\| \le
\frac{24\sqrt{6\ln N}}{5\pi e^2 N}(1+o(1)) \left\|f \right\|,\ \forall f\in L^2(\mr).
\end{equation}
\et

\bp The idea of the proof is to apply Lemma \ref{pea} (ii) on the dual frames
\bes \{f_k\}_{k\in I}:= \{E_{mb}T_{na}\, g_N\}_{m,n\in\mz}, \, \{g_k\}_{k\in I}:=
\{E_{mb}T_{na}\, h_N\}_{m,n\in\mz}.\ens Letting
$\{\widetilde f_k\}_{k\in I}:= \{E_{mb}T_{na}\, \widetilde g\}_{m,n\in\mz},$ the condition
\eqref{60812a} is satisfied with $\mu^2 = \widehat{B}(p_N, a,b),$ where $p_N$ is the function
defined in \eqref{60802d} and $\widehat{B}(p_N, a,b)$ is given by \eqref{16g-0}.
Thus, the term $\mu \sqrt{B_g}$ in \eqref{70220a} corresponds
precisely to
$$\left(\widehat{B}(p_N, a,b)B(h_N, a,b)\right)^{1/2}.$$
In order to estimate $\widehat{B}(p_N, a,b)$
we will use some of the calculations from the proof of Lemma \ref{60802c}.
Fix any integer $N\geq 14$.
Applying \eqref{16g-21} with $a:=\sqrt{\frac{12}{N}}$ yields
that
\begin{eqnarray*}
P_N(a)&=& P_N\left(\sqrt{\frac{12}{N}}\right) \\
&=&
\left(1+2 \sqrt{\frac{12}{N}}\frac{\sqrt{\ln N}}{\pi} \right)\frac{4}{5e^{2}N}\left(1+\frac{17\ln N}{7N}\right)\\
&&+
\left(2+2 \sqrt{\frac{12}{N}}\left( \frac{\sqrt{N}}{4\sqrt{3}}-\frac{\sqrt{\ln N}}{\pi} \right)\right)
\frac{4(\ln N)^2}{5N^3}\left(1+\frac{17\ln N}{7N}\right)\\
&&
+2\left( \frac{e^{-\frac{\pi^2 N}{24}}}{1-e^{-\frac{\pi^2 N}{6}}}  +
\left(\frac{2}{\pi}\right)^N\left(1+\frac{1}{2(N-1)} \right)\right)\\
&=& \frac{16\sqrt{3}\sqrt{\ln N}}{5\pi e^2 N\sqrt{N}}(1+o(1)).
\end{eqnarray*}
Similarly, we have
\begin{eqnarray*}
P_N(1/b)= P_N\left((2N-1)\sqrt{\frac{12}{N}}\right)
=\frac{32\sqrt{3}\sqrt{\ln N}}{5\pi e^2\sqrt{N}}(1+o(1)).
\end{eqnarray*}
Inserting this into \eqref{16g-22} yields
\begin{eqnarray*}
\widehat{B}(p_N,a,b)\leq \frac{P_N(a) P_N(1/b)}{a}=\sqrt{\frac{N}{12}}P(a) P(1/b)
= \frac{2^8\sqrt{3}\ln N}{5^2\pi^2 e^4   N\sqrt{N}}(1+o(1)).
\end{eqnarray*}
This together with \eqref{16g-32} implies that
$$\widehat{B}(p_N,a,b) B(h_N,a,b) \leq \frac{2^7 3^3 \ln N}{5^2 \pi^2 e^4 N^2}(1+o(1)).$$
Now \eqref{16g-34} is a consequence of Lemma \ref{pea} (ii).\ep

Theorem \ref{60812f} has a clear potential for applications, as we will explain in detail now.  Indeed, the known (exact) dual windows
associated with the Gabor frame generated by the Gaussian have a very complicated structure.
To the best of our knowledge, the only explicit expressions for dual windows associated with
the Gaussian were obtained by Janssen \cite{J2,Jan9}.  Consider for the moment the Gaussian
on the form $g(x)=2^{1/4}e^{-\pi x^2},$ and fix
$a,b>0$ such that $ab<1.$ Furthermore, let
\bes K:= \sum_{k\in \mz} (-1)^k (2k+1)e^{-\pi a(k+1/2)^2}, \quad
\mbox{erfc}(x):= \frac{2}{\sqrt{\pi}}\int_x^\infty e^{-s^2} \,ds.\ens
Then Janssen showed in \cite{J2} that for any $\epsilon>0, \ \epsilon < 1-ab,$ the function
\bee \label{15} \widetilde{g}_\epsilon(x):= 2^{-1/4} bK^{-1} e^{\pi x^2} \sum_{k\in \mz}
(-1)^k e^{-\pi a (k+1/2)^2/b} \ \mbox{erfc} [(x-(k+1/2)a)\sqrt{\pi/\epsilon}]\ene
generates a dual frame of $\mts$.
An explicit dual window is calculated in \cite{Jan9} for the case where $(ab)^{-1}$ is an even integer (the calculations are more involved when
$(ab)^{-1}$ is an odd integer).
On the other hand, as already described the functions $h_N$ in \eqref{16g-31} have an extremely simple and explicit structure.
Theorem \ref{60812f} shows that for certain sufficiently small
values of the product $ab$ we can get arbitrary close to perfect
reconstruction by replacing the exact dual window by the function $h_N$ for a sufficiently
large value of $N\in \mn.$  For example, doing so with a tolerance that is smaller then the ``machine precision"
(in a slightly informal language) will be as good as having perfect reconstruction.

Note that in Theorem \ref{60812f} each choice of $N\in \mn$ leads to new and fixed values
for the translation parameter $a$ and the modulation parameter $b,$ given by typically irrational
numbers.  For general frames this can be problematic, due to inevitable roundoff errors that
will occur in practice.  However, the Gaussian and the (scaled) B-splines $B_N, \, N\ge 2,$ belong
to the Feichtinger algebra ${\cal S}_0,$ which implies that the Bessel bound for the
Gabor system $\{E_{mb}T_{na}h_N\}_{m,n\in \mz}$ depends continuously on the parameters
$a,b>0,$ cf. the results by Feichtinger and Kaiblinger \cite{FK}.  Thus, Theorem \ref{60812f} is less sensitive towards roundoff errors than one could
think at a first glance.

\section{Convergence in $L^p$-spaces} \label{60805a}

In this section we derive some easy consequences of our approach in Section \ref{060801a}.
In Theorem \ref{60817a} we prove that  the considered scalings of the B-splines converge to the
Gaussian in $L^p(\mr)$ for $p\in[1,\infty];$ this was also obtained by Goh et al. in
\cite{GGL}, while  Unser et al. \cite{AU} proved it for $p\in [2, \infty].$
First we prove that the Fourier transform of the scaled
B-splines converge to the Gaussian in $L^q(\mr)$ for $q\in [1,
\infty];$ the range $q\in [1, \infty[$ was first considered by Unser et al. \cite{AU}.

\bt \label{16g-37}
Let $q\in[1,\infty].$
Then
$$
\mathcal{F}\left(\sqrt{\frac{N}{12}} \, B_N\left(\sqrt{\frac{N}{12}}\, \cdot\right)\right)(\ga)
=\left(\frac{\sin\left(\pi\sqrt{\frac{12}{N}}\, \ga \right)}{\pi\sqrt{\frac{12}{N}}\,\ga}\right)^{N}
\rightarrow
e^{-2\pi^2 \ga^2}$$
in $L^q(\mr)$ as $N\rightarrow \infty$. \et
\bp
Fix any integer $N\geq 14$.
We will distinguish between the cases $q=\infty$ and  $q\in [1,\infty[$.

(1) The case $q=\infty$:
By \eqref{16g-8}, we see that
$$\left| e^{-2\pi^2 \ga^2}
-\left(\frac{\sin\left(\pi\sqrt{\frac{12}{N}}\, \ga \right)}{\pi\sqrt{\frac{12}{N}}\,\ga}\right)^{N}\right|
\leq e^{-\frac{\pi^2 N}{24}}
+\left(\frac{2}{\pi}\right)^N, \ |\ga|\geq \frac{\sqrt{N}}{4\sqrt{3}}.$$
This together with \eqref{16g-7} implies that
\begin{equation*}
\sup_{\ga\in\mr}\left| e^{-2\pi^2 \ga^2}
-\left(\frac{\sin\left(\pi\sqrt{\frac{12}{N}}\, \ga \right)}{\pi\sqrt{\frac{12}{N}}\,\ga}\right)^{N}\right|
\leq \max\left\{\frac{4}{5e^{2}N}\left(1+\frac{17\ln N}{7N}\right), e^{-\frac{\pi^2 N}{24}}
+\left(\frac{2}{\pi}\right)^N \right\} \rightarrow 0
\end{equation*}
as $N\rightarrow\infty.$

(2) The case $q\in [1,\infty[$:
By
\eqref{16g-7} and \eqref{16g-8}, we have
\begin{eqnarray*}
&& \int_{-\infty}^\infty
\left| e^{-2\pi^2 \ga^2}
-\left(\frac{\sin\left(\pi\sqrt{\frac{12}{N}}\, \ga \right)}
{\pi\sqrt{\frac{12}{N}}\,\ga}\right)^{N}\right|^q d\ga\\
&&\leq
\int_{|\ga|< \frac{\sqrt{\ln N}}{\pi}}
\left(\frac{4}{5e^{2}N}\left(1+\frac{17\ln N}{7N}\right)\right)^q d\ga+
\int_{\frac{\sqrt{\ln N}}{\pi}\leq |\ga|\leq \frac{\sqrt{N}}{4\sqrt{3}}}
\left(\frac{4(\ln N)^2}{5N^3}\left(1+\frac{17\ln N}{7N}\right)\right)^q d\ga \\
&&+
\int_{|\ga|\geq \frac{\sqrt{N}}{4\sqrt{3}}}
\left(e^{-2\pi^2 \ga^2} +
\left|\frac{1}{\pi\sqrt{\frac{12}{N}}\,\ga}\right|^{N}\right)^q d\ga.
\end{eqnarray*}
Using that $(a+b)^q\leq \max\{(2a)^q, (2b)^q  \}\leq 2^q(a^q+b^q)$,
we have
\begin{eqnarray*}
&&  \int_{|\ga|\geq \frac{\sqrt{N}}{4\sqrt{3}}}
\left(e^{-2\pi^2 \ga^2} +
\left|\frac{1}{\pi\sqrt{\frac{12}{N}}\,\ga}\right|^{N}\right)^q d\ga \\
&&\leq 2^q \left(\int_{|\ga|\geq \frac{\sqrt{N}}{4\sqrt{3}}}
e^{-2\pi^2 \ga^2 q}d\ga+
 \int_{|\ga|\geq \frac{\sqrt{N}}{4\sqrt{3}}}
\left|\frac{1}{\pi\sqrt{\frac{12}{N}}\,\ga}\right|^{Nq}d\ga\right).
\end{eqnarray*}
Direct calculations show that
\begin{eqnarray*}
&& \int_{|\ga|\geq \frac{\sqrt{N}}{4\sqrt{3}}}e^{-2\pi^2 \ga^2 q}d\ga
\leq 2\int_{\frac{\sqrt{N}}{4\sqrt{3}}}^\infty e^{-2\pi^2 \ga q}d\ga     =
\frac{e^{-\left(2\pi^2 q \frac{\sqrt{N}}{4\sqrt{3}}\right)}}{\pi^2 q} \
\left(\text{if} \ \frac{\sqrt{N}}{4\sqrt{3}}\geq 1\right);\\
&&\int_{|\ga|\geq \frac{\sqrt{N}}{4\sqrt{3}}}
\left|\frac{1}{\pi\sqrt{\frac{12}{N}}\,\ga}\right|^{N q}d\ga
=\left( \frac{2}{\pi} \right)^{Nq} \frac{\sqrt{N}}{2\sqrt{3}(N q-1)} \
\left(\text{if} \ Nq>1\right).
\end{eqnarray*}
Hence
\begin{eqnarray}
&& \int_{-\infty}^\infty
\left| e^{-2\pi^2 \ga^2}
-\left(\frac{\sin\left(\pi\sqrt{\frac{12}{N}}\, \ga \right)}
{\pi\sqrt{\frac{12}{N}}\,\ga}\right)^{N}\right|^q d\ga \nonumber \\
 && \leq
2\left( \left(\frac{4}{5e^{2}N}  \right)^q \frac{\sqrt{\ln N}}{\pi}  +
 \left(\frac{4(\ln N)^2}{5N^3}  \right)^q
\left(\frac{\sqrt{N}}{4\sqrt{3}}-\frac{\sqrt{\ln N}}{\pi} \right) \right)\left(1+\frac{17\ln N}{7N}\right)^q
\nonumber\\
&&+2^q\left( \frac{ e^{-\left(2\pi^2 q \frac{\sqrt{N}}{4\sqrt{3}}\right)}}{\pi^2 q}
+ \left( \frac{2}{\pi} \right)^{N q} \frac{\sqrt{N}}{2\sqrt{3}(N q-1)}\right) \nonumber \\
&& \leq C_q\frac{\sqrt{\ln N}}{N^q}, \label{16g-24}
\end{eqnarray}
for some  $C_q>0$, and therefore tends to zero as $N\rightarrow\infty$.
\ep

We are now ready to prove that the considered scaled versions of the B-splines $B_N$ converge to
the Gaussian in $L^p(\mr), \, p\in [1, \infty]$ as $N\to \infty.$

\bt \label{60817a}
Let $p\in [1,\infty].$
Then
$$
\sqrt{\frac{N}{12}} \, B_N\left(\sqrt{\frac{N}{12}}\, x \right)
\rightarrow \frac{1}{\sqrt{2\pi}}e^{-\frac{x^2}{2}} $$
in $L^p(\mr)$ as $N\rightarrow \infty$. \et
\bp
We will distinguish between the cases  $p=\infty$ and  $p\in [1,\infty[$.

(1) The case $p=\infty$: Using the function $p_N$ in \eqref{60802d} and
its Fourier transform in \eqref{16g-12}, the inequality
$\| p_N \|_{\infty} \leq \|\widehat p_N \|_{L^1(\mr)}$
and \eqref{16g-24} with $q=1$ imply that
\begin{eqnarray}
\left\|
\frac{1}{\sqrt{2\pi}}\, e^{-\frac{x^2}{2}}-
 \sqrt{\frac{N}{12}} \, B_N\left(\sqrt{\frac{N}{12}}\, x \right) \right\|_{\infty}
&\leq& \int_{-\infty}^\infty \left| e^{-2\pi^2 \ga^2}
-\left(\frac{\sin\left(\pi\sqrt{\frac{12}{N}}\, \ga \right)}
{\pi\sqrt{\frac{12}{N}}\,\ga}\right)^{N}\right| d\ga \nonumber \\
&\leq&  C_1 \frac{\sqrt{\ln N}}{N}, \label{16g-25}
\end{eqnarray}
which tends to zero as $N\rightarrow\infty$.

(2) The case $p\in[1,\infty[$: Since $\supp \, B_N(\sqrt{\frac{N}{12}}\, \cdot)=[-\sqrt{3N},\sqrt{3N}]$,
we see that
\begin{eqnarray*}
&& \int_{-\infty}^\infty
\left|\frac{1}{\sqrt{2\pi}}\, e^{-\frac{x^2}{2}}-\sqrt{\frac{N}{12}}
 \, B_N\left(\sqrt{\frac{N}{12}}\, x \right)
 \right|^p dx\\
&&\leq
\int_{|x|\leq \sqrt{3N}} \left|
\frac{1}{\sqrt{2\pi}}\, e^{-\frac{x^2}{2}}-\sqrt{\frac{N}{12}}
 \, B_N\left(\sqrt{\frac{N}{12}}\, x \right) \right|^p dx+
\int_{|x|\geq \sqrt{3N}} \left|\frac{1}{\sqrt{2\pi}}\, e^{-\frac{x^2}{2}} \right|^p dx.
\end{eqnarray*}
By \eqref{16g-25}, we have
$$\int_{|x|\leq \sqrt{3N}} \left|\frac{1}{\sqrt{2\pi}}\, e^{-\frac{x^2}{2}}-\sqrt{\frac{N}{12}}
 \, B_N\left(\sqrt{\frac{N}{12}}\, x \right) \right|^p dx
\leq  2\sqrt{3N}
\left(\frac{C_1\sqrt{\ln N}}{N}\right)^p =
\frac{ 2C_1^p \sqrt{3} \left(\sqrt{\ln N}\right)^p }{N^{p-\frac{1}{2}}}.
$$
Also direct calculations show
$$
\int_{|x|\geq \sqrt{3N}} \left|\frac{1}{\sqrt{2\pi}}\, e^{-\frac{x^2}{2}} \right|^p dx
\leq
\int_{|x|\geq \sqrt{3N}} \left|\frac{1}{\sqrt{2\pi}}\, e^{-\frac{x}{2}} \right|^p dx =
\frac{4}{p} \left(\frac{1}{\sqrt{2\pi}}\right)^p e^{-\frac{p\sqrt{3N}}{2}}.
$$
Hence
\begin{eqnarray*}
\int_{-\infty}^\infty
\left|\sqrt{\frac{N}{12}} \, B_N\left(\sqrt{\frac{N}{12}}\, x \right)
-\frac{1}{\sqrt{2\pi}}\, e^{-\frac{x^2}{2}} \right|^p dx
\leq \frac{ 2C_1^p \sqrt{3} \left(\sqrt{\ln N}\right)^p }{N^{p-\frac{1}{2}}}
+\frac{4}{p} \left(\frac{1}{\sqrt{2\pi}}\right)^p e^{-\frac{p\sqrt{3N}}{2}}
\rightarrow 0
\end{eqnarray*}
 as $N\rightarrow\infty.$
\ep

Consider again the function $g_N$ in \eqref{60802g}.
As a consequence of our results we obtain that a certain scaled version of the Fourier transform
of $g_N$ ``behave very similar as $g_N$" as $N\to \infty.$  Recall that for $a>0$ the
scaling operator $D_a$ on $\ltr$ is defined by
$D_af(x)= a^{-1/2}f(a^{-1}x).$ The scaling operator is unitary.

\bc
Let $p\in [1,\infty].$
Then
$
g_N- D_{2\pi}
\mathcal{F}g_N
\rightarrow 0
$
in $L^p(\mr)$ as $N\rightarrow \infty$.
\ec
The result follows immediately by letting $\ga=\frac{x}{2\pi}$ in Theorem \ref{16g-37}.


\vspace{1cm}\noindent{\bf Acknowledgments:} The authors would like to express their gratitude to Guido Janssen for many
detailed comments on an earlier version of the paper.  The comments improved
the presentation and also lead to much better estimates for some of the
Bessel bounds. Furthermore O. C.
thanks Say Song Goh for making him aware of the papers
\cite{Bricks,GGL,AU}. The authors  thank NIMS for
support and hospitality during their visit in February 2017.
This research was supported by Basic Science Research Program
through the National Research Foundation of Korea(NRF) funded by
the Ministry of Education(2016R1D1A1B02009954) and also by
National Institute for Mathematical Sciences (NIMS) (A23100000).

{\bf \noindent Ole Christensen\\
Department of Applied Mathematics and Computer Science\\
Technical University of Denmark,
Building 303,
2800 Lyngby  \\
Denmark, Email: ochr@dtu.dk

\vspace{.1in}\noindent Hong Oh Kim \\
Department of Mathematical Sciences, UNIST\\
50 UNIST-gil, Ulsan 44919\\
Republic of Korea, Email: hkim2031@unist.ac.kr

\vspace{.1in} \noindent Rae Young Kim \\
Department of Mathematics, Yeungnam University\\
280 Daehak-Ro, Gyeongsan, Gyeongbuk 38541\\
Republic of Korea, Email:  rykim@ynu.ac.kr }


\begin{thebibliography}{10}

\frenchspacing


\bibitem{Bricks} Bricks, R.:  On the convergence of derivatives of B-splines.
    Comput. Appl. Math. {\bf 27}(1), 80--92 (2008)



\bibitem{CC3} Casazza, P. G., Christensen, O. and Janssen, A.J.E.M.:
{\it Weyl-Heisenberg frames, translation invariant systems, and
the Walnut representation.} J. Funct. Anal. {\bf 180} (2001),
85--147.

\bibitem{C19} Christensen, O.: Pairs of dual Gabor frames with compact support and desired
frequency localization. Appl. Comput. Harmon. Anal. {\bf 20},
403--410 (2006)


\bibitem{CB} Christensen, O.:  An introduction to frames and Riesz bases. Second expanded
edition. Birkh\"auser (2016)


\bibitem{CH} Christensen, O. and Heil, C.: {\it Perturbations of Banach frames
and atomic decompositions.} Math. Nach. {\bf 185} (1997),  33--47.

\bibitem{CR} Christensen, O., Kim, R.Y.:  On dual Gabor
frame pairs generated by polynomials. J. Fourier Anal. Appl. {\bf 16}, 1--16 (2010)





\bibitem{CLau} Christensen, O., Laugesen, R.:  Approximately dual frames in
Hilbert spaces and applications to Gabor frames. Sampl. Theory
Signal Image Process. {\bf 9}, 77--90 (2011)


\bibitem{FK} Feichtinger, H.G., Kaiblinger, N.:  Varying
the time-frequency lattice of Gabor frames.  Trans. Amer. Math. Soc. {\bf 356},   2001--2023 (2004)

\bibitem{GGL} Goh, S.S., Goodman, T.N.T.,  Lee, S.L.:
 Appell sequences, continuous wavelet transforms, and series expansions.
 Appl. Comput. Harmon. Anal. to appear 2016



\bibitem{J2} Janssen, A.J.E.M.:  Signal analytical proofs of two basic
results on lattice expansions. Appl. Comput. Harmon. Anal. {\bf 1},
350-354 (1994)


\bibitem{Jan9} Janssen, A.J.E.M.: {\it Some Weyl-Heisenberg frame bound calculations.}
Indag. Math {\bf 7} (1996), 165--183.


\bibitem{LemNiel} Lemvig, J.,  Nielsen, K.H.:   Counterexamples to the B-spline
conjecture for Gabor frames. J. Fourier Anal. Appl. {\bf 22}, 1440-1451 (2016)






\bibitem{Ly} Lyubarskii, Y.:  Frames in the Bargmann space of
entire functions. Adv. Soviet Math. {\bf 11}, 167--180 (1992)






\bibitem{Se2} Seip, K.: Density theorems for sampling and interpolation in the Bargmann-Fock
space I. J. Reine Angew. Math. {\bf 429},  91--106 (1992)


\bibitem{SW} Seip, K., Wallsten, R.:  Density theorems for sampling and interpolation in the
Bargmann-Fock space II. J. Reine Angew. Math. {\bf 429},
107--113 (1992)


\bibitem{AU} Unser, M., Aldroubi, A., Eden, M.:  On the asymptotic convergence
of B-spline wavelets to Gabor functions. IEEE Trans. Inform. Theory {\bf 38}(2), 864--872 (1992)

\end{thebibliography}
\end{document}